\documentclass[12pt,leqno]{amsproc}
\usepackage{amsmath,amssymb}
\title{Generalized Divisors and Biliaison}
\author{Robin Hartshorne}
\date{}
\address{Department of Mathematics\\
University of California\\
Berkeley, CA\ \ 94720-3840}
\email{robin@math.berkeley.edu}

\begin{document}

\begin{abstract}
We extend the theory of generalized divisors so as to work on any scheme $X$
satisfying the condition $S_2$ of Serre.  We define a generalized notion of
Gorenstein biliaison for schemes in projective space.  With this we give a new
proof in a stronger form of the theorem of Gaeta, that standard determinantal schemes are in the
Gorenstein biliaison class of a complete intersection.

We also show, for schemes of codimension three in ${\mathbb P}^n$, that the relation of Gorenstein biliaison is equivalent to the relation of even strict Gorenstein liaison.
\end{abstract}

\maketitle

\setcounter{section}{-1}
\section{Introduction}
\label{sec0}

In this paper we generalize further the theory of generalized divisors
introduced in \cite{GD} by partially removing the Gorenstein hypotheses.  This,
we feel, puts the theory in its natural state of generality.  The main
difference is that instead of requiring the sheaves of ideals defining a
generalized divisor to be reflexive, we require only the condition $S_2$ of
Serre.  If a scheme $X$ satisfies $G_1$ and $S_1$, then a coherent sheaf is
reflexive if and only if it satisfies $S_2$ \cite[1.9]{GD}.  Here we show that
if $X$ satisfies $S_1$ only, then a coherent sheaf satisfies $S_2$ if and only
if it is $\omega$-{\em reflexive}:  this means that the natural map ${\mathcal
F} \rightarrow {\mathcal H}om({\mathcal H}om({\mathcal F},\omega),\omega)$ is
an isomorphism, where $\omega$ is the canonical sheaf.  With this weaker
condition we are able to establish a theory of generalized divisors on schemes
$X$ satisfying only the condition $S_2$.

We apply this theory to define a notion of generalized biliaison for schemes in
projective space.  Let $D$ be a (generalized) divisor on an ACM scheme $X$ in
${\mathbb P}_k^n$.  If $D' \sim D + mH$, meaning $D'$ is linearly equivalent to
$D$ plus $m$ times the hyperplane section $H$ of $X$, we say $D'$ is obtained
by an {\em elementary biliaison} from $D$.  We call {\em biliaison} the
equivalence relation generated by the elementary biliaisons.

If we do biliaisons using only complete intersection schemes $X$ in ${\mathbb
P}^n$, the resulting notion of biliaisons is equivalent to even complete
intersection liaison (CI-liaison) \cite[4.4]{GD}.  If we do biliaisons using an
ACM scheme $X$ satisfying $G_0$, we will show $(3.5)$ that any such
biliaison (called $G$-biliaison) is an even Gorenstein liaison.  We do not know
if the converse is true.  If we use arbitrary ACM schemes $X$, we obtain a
notion of biliaison that is possibly more general than $G$-biliaison.
Note that even this more general type of biliaison preserves the Rao modules up
to shift $(3.2)$.

As an application we give a new proof $(4.1)$ of the theorem of KMMNP--Gaeta \cite[Sec.
3]{KMMNP} using biliaisons.  I would like to thank Marta Casanellas for
explaining the old proof and helping to discover the new proof given here.

Examining the proof of $(3.5)$ we see that the Gorenstein linkages there are all of a special kind: they use only arithmetically Gorenstein schemes of the form $M+mH$ on some ACM scheme satisfying $G_0$ (cf.~$3.3$).  These we call {\em strict} Gorenstein linkages, and so $(3.5)$ actually tells us that every $G$-biliaison is an even strict Gorenstein liaison.  In section~\ref{sec5} we prove a partial converse: Every even strict Gorenstein liaison of codimension~$3$ subschemes of ${\mathbb P}^n$ is a Gorenstein biliaison $(5.1)$.

\section{$\omega$-Reflexive modules}
\label{sec1}

We will need some well-known results about the canonical module, or dualizing
module as it is sometimes called, of a ring or scheme.  We restrict our
attention to equidimensional {\em embeddable} noetherian rings and schemes.  For
a ring $A$, this means that it is a quotient of a regular ring.  For a scheme
$X$, it means that it can be embedded as a closed subscheme of a regular scheme. 
This includes all quasiprojective schemes over a field, which will be our most
common application.

An equidimensional embeddable ring or scheme always has a canonical
module or sheaf unique up to isomorphism.  It is finitely generated (resp.\
coherent).  Its formation commutes with localization, and with completion of a
local ring.  If the ring $A$ is a quotient of a regular ring $P$, and $r$ is
the difference of dimensions, then the canonical module $\omega$ of $A$ can be
obtained as $\omega = \mbox{Ext}_P^r(A,P)$, and similarly for a closed
subscheme $X$ of a regular scheme $P$.  If $A$ is a Cohen--Macaulay ring, then
$\omega$ is a Cohen--Macaulay module of the same dimension as $A$, and for any
maximal Cohen--Macaulay module $M$, the natural map $M \rightarrow
\mbox{Hom}_A(\mbox{Hom}_A(M,\omega),\omega)$ is an isomorphism.  For
references see \cite{HK} for the case of Cohen--Macaulay rings; \cite[II.7]{AG}
for the case of projective schemes; and see also \cite{A}.

We will expand these results somewhat by weakening their hypotheses to suit our
situation.  We define a module $M$ over a ring $A$ (as above) to be
$\omega$-{\em reflexive} if the natural map $M \rightarrow
\mbox{Hom}_A(\mbox{Hom}_A(M,\omega),\omega)$ is an isomorphism.  Sometimes we
will denote by
$M^{\omega}$ the module $\mbox{Hom}_A(M,\omega)$, and call it the
$\omega$-{\em dual} of $M$.

\bigskip
\noindent
{\bf Lemma 1.1.} {\em If $A$ is a local ring of dimension $0$, every finitely
generated module $M$ is $\omega$-reflexive.}

\bigskip
\noindent
{\bf {\em Proof.}} Since $A$ is Cohen--Macaulay, this follows from \cite[6.1]{HK}. 
It also follows from the local duality theorem, which says in this case that
$M^{\omega}$ is the dual of $H_{\mathfrak m}^0(M) = M$, so that
$M^{\omega\omega}$ is the double dual, which is isomorphic to $M$.

\bigskip
\noindent
{\bf Lemma 1.2.} {\em For any local ring $A$, the module $\omega_A$ satisfies the
condition $S_1$ of Serre.}

\bigskip
\noindent
{\bf {\em Proof.}} Write $A$ as a quotient of a regular local ring $P$ of codimension
$r$.  Then $\omega_A = \mbox{Ext}_P^r(A,P)$.  By reason of dimension and local
duality on $P$, the functor $\mbox{Ext}_P^r(\cdot,P)$ is contravariant and
left-exact on $A$ modules.  If $\dim A = 0$, there is nothing to prove.  If
$\dim A \ge 1$, let 
$x
\in {\mathfrak m}_A$ be an element such that
$\dim A/xA < \dim A$.  Then from the sequence
\[
A \stackrel{x}{\rightarrow} A \rightarrow A/xA \rightarrow 0
\]
we obtain
\[
0 = \mbox{Ext}_P^r(A/xA,P) \rightarrow \omega_A \stackrel{x}{\rightarrow}
\omega_A.
\]
Thus $\omega_A$ has depth $1$.  Since formation of $\omega$ commutes with
localization, we conclude that $\omega_A$ satisfies $S_1$.

\bigskip
\noindent
{\bf Lemma 1.3.} {\em If a local ring $A$ satisfies $S_1$, then $\omega_A$
satisfies $S_2$.}

\bigskip
\noindent
{\bf {\em Proof.}} Write $A$ as a quotient of a regular local ring $P$ of codimension
$r$, as before.  Let $x \in {\mathfrak m}_A$ be a non-zero-divisor so that $B =
A/xA$ has dimension one less.  Then from the exact sequence
\[
0 \rightarrow A \stackrel{x}{\rightarrow} A \rightarrow B \rightarrow 0
\]
and $(1.2)$ we obtain
\[
0 \rightarrow \omega_A \stackrel{x}{\rightarrow} \omega_A \rightarrow
\mbox{Ext}_p^{r+1}(B,P) = \omega_B.
\]
Since $\omega_B$ satisfies $S_1$ by $(1.2)$, we see that if $\dim A \ge
2$, then $\omega_A$ has depth $\ge 2$.  Hence $\omega_A$ satisfies $S_2$.

\bigskip
\noindent
{\bf Lemma 1.4.} {\em Let $A$ be a one-dimensional local Cohen--Macaulay ring. 
Then a finitely generated module $M$ is $\omega$-reflexive if and only if it has
depth $1$.}

\bigskip
\noindent
{\bf {\em Proof.}} Since $\omega$ has depth $1$ $(1.2)$ so does the $\omega$-dual of
any module.  If $M$ is reflexive, it is the $\omega$-dual of $M^{\omega}$ and
so has depth $1$.  The converse is \cite[6.1]{HK}.

\bigskip
\noindent
{\bf Proposition 1.5.} {\em Let $A$ be a local ring satisfying $S_1$.  A
finitely generated module $M$ is $\omega$-reflexive if and only if it satisfies
$S_2$.}

\bigskip
\noindent
{\bf {\em Proof.}} First we show that the $\omega$-dual of any module $N$ will
satisfy $S_2$.  Write $N$ as a cokernel of a map of free modules
\[
L_1 \rightarrow L_0 \rightarrow N \rightarrow 0.
\]
Taking $\omega$-duals and the image of the second map, we obtain
\[
0 \rightarrow N^{\omega} \rightarrow L_0^{\omega} \rightarrow K \rightarrow 0
\]
where $K$ is a submodule of $L_1^{\omega}$.  Now $L_0^{\omega}$ and
$L_1^{\omega}$ are direct sums of copies of $\omega$, so satisfy $S_2$ by
$(1.3)$.  Hence $K$ satisfies $S_1$, and then from the exact sequence it
follows that $N^{\omega}$ satisfies $S_2$.  In particular, any
$\omega$-reflexive module satisfies $S_2$.

Conversely, suppose $M$ satisfies $S_2$.  The map $\alpha: M \rightarrow
M^{\omega\omega}$ is an isomorphism in codimension $0$, by $(1.1)$, so the
kernel of $\alpha$ must have support of codimension $\ge 1$.  Since $M$
satisfies $S_1$, there is no kernel.  Thus we can write
\[
0 \rightarrow M \stackrel{\alpha}{\rightarrow} M^{\omega\omega} \rightarrow R
\rightarrow 0.
\]
Now, since $\alpha$ is an isomorphism in codimension $1$, by $(1.4)$, the
module $R$ must have support of codimension $\ge 2$.  Since both $M$ and
$M^{\omega\omega}$ satisfy $S_2$, this is impossible (cf.\ proof of
\cite[1.9]{GD}) so $R = 0$ and $\alpha$ is an isomorphism.

\bigskip
\noindent
{\bf Corollary 1.6.} {\em Let $A$ satisfy $S_1$.  The $\omega$-dual of any
module is $\omega$-reflexive.}

\bigskip
\noindent
{\bf {\em Proof.}} This follows from the first step of the proof of $(1.5)$.

\bigskip
\noindent
{\bf Corollary 1.7.} {\em If $A$ itself satisfies $S_2$, then the natural map
$A \rightarrow \mbox{\em Hom}_A(\omega,\omega)$ is an isomorphism.}

\bigskip
\noindent
{\bf Remark 1.8.} Using the same arguments as in \cite[1.11,1.12]{GD} we see
that if ${\mathcal F}$ is a coherent sheaf satisfying $S_2$ on a scheme $X$
satisfying $S_1$, then ${\mathcal F}$ is {\em normal} in the sense of Barth
\cite[1.6]{SRS}, namely for any open set $U$ and any closed subset $Y \subseteq
U$ of codimension $\ge 2$, the restriction map ${\mathcal F}(U) \rightarrow
{\mathcal F}(U-Y)$ is bijective.  In fact this condition characterizes $S_2$,
if we assume $S_1$.

If ${\mathcal F}$ is a coherent sheaf satisfying $S_1$ only, then it is easy to
see that the set $Y$ of points of $X$ where it does not satisfy $S_2$ is a
closed subset of codimension $\ge 2$, and that the double $\omega$-dual
${\mathcal F}^{\omega\omega}$ can be identified with $j_*({\mathcal F}|_{X-Y})$
where $j: X - Y \rightarrow X$ is the inclusion.  Thus the double $\omega$-dual
can be regarded as the $S_2$-ification of the sheaf.

It also follows naturally that for $Y \subseteq X$ closed of codimension $\ge
2$, the category of coherent sheaves satisfying $S_2$ on $X$ is equivalent by
restriction to the analogous category on $X-Y$.

\bigskip
\noindent
{\bf Remark 1.9.} To see the connection between the properties reflexive and
$\omega$-reflexive, note that the proof of \cite[1.9]{GD} shows that a
reflexive module over a ring $A$ satisfying $S_2$ also satisfies $S_2$.  So we
see that if $A$ satisfies $S_2$, then a reflexive module is also
$\omega$-reflexive.  The converse is not true without the $G_1$ hypothesis. 
For example, if $X$ is the union of the three coordinate axes in ${\mathbb
A}^3$, a scheme that satisfies $G_0$ but not $G_1$, the canonical sheaf
$\omega$ is $\omega$-reflexive by $(1.4)$, but is easily seen not to be
reflexive.  On the other hand, the proof of \cite[1.9]{GD} does show that if
$X$ satisfies $S_2$, and ${\mathcal F}$ satisfies $S_2$ and is reflexive in
codimension $\le 1$, then ${\mathcal F}$ is reflexive.

\section{Generalized divisors}
\label{sec2}

Let $X$ be a noetherian, equidimensional, embeddable scheme satisfying the
condition $S_2$ of Serre.  We develop the theory of generalized divisors as in
\cite[\S 2]{GD}, noting the differences in our more general setting.

Let ${\mathcal K}_X$ be the sheaf of total quotient rings on $X$
\cite[2.1]{GD}.  A {\em fractional ideal} is a subsheaf ${\mathcal I} \subseteq
{\mathcal K}_X$ that is a coherent sheaf of ${\mathcal O}_X$-modules.  It is
{\em nondegenerate} if for each generic point $\eta \in X$, ${\mathcal
I}_{\eta} = {\mathcal K}_{X,\eta}$.

\bigskip
\noindent
{\bf Definition.}  Let $X$ be a scheme (as above) satisfying $S_2$.  A {\em
generalized divisor} on $X$ is a nondegenerate fractional ideal ${\mathcal I}$
satisfying the condition $S_2$ as a sheaf of ${\mathcal O}_X$-modules.  It is
{\em effective} if ${\mathcal I} \subseteq {\mathcal O}_X$.  We say the
generalized divisor ${\mathcal I}$ is {\em principal} if ${\mathcal I} = (f)$
for some global section $f \in {\mathcal K}_X$.  We say it is {\em Cartier} if
${\mathcal I}$ is an invertible ${\mathcal O}_X$-module.  We say it is {\em
almost Cartier} if there exists a closed subset $Z \subseteq X$ of codimension
$\ge 2$ so that ${\mathcal I}|_{X-Z}$ is Cartier.  We say it is {\em reflexive}
if ${\mathcal I}$ is a reflexive ${\mathcal O}_X$-module.

Note that ${\mathcal I}$ is Cartier if and only if it is locally principal
\cite[2.3]{GD}.  Note that an almost Cartier divisor is reflexive $(1.9)$ and
that the sheaf ${\mathcal I}$ being reflexive implies the condition $S_2$.

\bigskip
\noindent
{\bf Proposition 2.1.} {\em With $X$ satisfying $S_2$, as above, the effective
generalized divisors are in one-to-one correspondence with closed subschemes $Y
\subseteq X$ of pure codimension one with no embedded points.}

\bigskip
\noindent
{\bf {\em Proof.}} Let $Y$ be a closed subscheme of $X$, defined by a sheaf of ideals
${\mathcal I}$, so that we have an exact sequence
\[
0 \rightarrow {\mathcal I} \rightarrow {\mathcal O}_X \rightarrow {\mathcal
O}_Y \rightarrow 0.
\]
To say that ${\mathcal I}$ is nondegenerate is equivalent to saying $Y$ has
codimension $\ge 1$.  Since $X$ satisfies $S_2$, to say that ${\mathcal I}$
satisfies $S_2$ is equivalent to saying that every associated prime of $Y$ has
codimension $1$ (cf.\ \cite[1.10]{GD}), i.e., that $Y$ is of pure codimension
$1$ with no embedded points.

\bigskip
\noindent
{\bf Definition.} For any coherent sheaf ${\mathcal F}$ of ${\mathcal
O}_X$-modules, let us denote ${\mathcal F}^{\sim}$ its double $\omega$-dual, so
that ${\mathcal F}^{\sim}$ satisfies $S_2$, where $\omega$ is the canonical
sheaf on $X$ $(1.8)$.  If ${\mathcal I} \subseteq {\mathcal K}_X$ is a
fractional ideal, then naturally ${\mathcal I}^{\sim}$ is also a fractional
ideal, and will satisfy $S_2$.  We may often denote a generalized divisor
${\mathcal I}$ by a letter $D$, and call ${\mathcal I}$ the ideal of $D$. 
Given two (generalized) divisors $D_1$ and $D_2$, with corresponding ideal
sheaves ${\mathcal I}_1,{\mathcal I}_2$, we define the {\em sum} $D_1 + D_2$ by
the fractional ideal $({\mathcal I}_1 \cdot {\mathcal I}_2)^{\sim}$.  We define
the negative $-D$ by $({\mathcal I}^{-1})^{\sim}$, where ${\mathcal I}^{-1} =
{\mathcal H}om({\mathcal I},{\mathcal O}_X)$.  We denote the divisor with ideal
${\mathcal I} = {\mathcal O}_X$ by $0$.

\bigskip
\noindent
{\bf Proposition 2.2.} {\em Let $X$ satisfy $S_2$.}
\begin{itemize}
\item[(a)] {\em Addition of divisors is associative and commutative.}
\item[(b)] {\em $D + 0 = D$ for all $D$.}
\item[(c)] {\em $-(-D) = D$ if and only if $D$ is reflexive.}
\item[(d)] {\em $D + (-D) = 0$ if and only if $D$ is almost Cartier.}
\item[(e)] {\em If $D$ is any divisor, and $E$ is almost Cartier, then $-(D+E)
= (-D) + (-E)$.}
\end{itemize}

\bigskip
\noindent
{\bf {\em Proof.}} (a) and (b) are obvious.  (c) follows from the fact ${\mathcal
I}^{-1} \cong {\mathcal I}^{\vee}$ as ${\mathcal O}_X$-modules \cite[2.2]{GD}. 
For (d) we follow the proof of \cite[2.5]{GD}, noting that at a point of
codimension $1$, every ideal is $\omega$-reflexive $(1.4)$, so that the
condition says ${\mathcal I} \cdot {\mathcal I}^{-1} = {\mathcal O}_X$,
which implies ${\mathcal I}$ reflexive there \cite[2.3]{GD}.  For (e) it is the
same proof as \cite[2.5]{GD}.

\bigskip
\noindent
{\bf Corollary 2.3.} {\em The set of almost Cartier divisors forms a group,
containing the subgroups of Cartier divisors and of principal divisors.  This
group acts on the set of all divisors.}

\bigskip
\noindent
{\bf Definition.} We say two divisors are {\em linearly equivalent} if one is
obtained from the other by adding a principal divisor.  We denote the
equivalence classes by the group $\mbox{Pic } X =$
Cartier divisors mod linear equivalence; the group $\mbox{APic } X =$ almost
Cartier divisors mod linear equivalence, and the set $\mbox{GPic } X =$
generalized divisors mod linear equivalence.

\bigskip
\noindent
{\bf Proposition 2.4.} {\em Two divisors $D_1$ and $D_2$ are linearly
equivalent if and only if their ideal sheaves ${\mathcal I}_1$ and ${\mathcal
I}_2$ are isomorphic as ${\mathcal O}_X$-modules.  Every coherent ${\mathcal
O}_X$-module that satisfies $S_2$ and is locally free of rank $1$ at every
generic point of $X$ is isomorphic to the ideal of some divisor.}

\bigskip
\noindent
{\bf {\em Proof.}} Indeed, an isomorphism $\varphi: {\mathcal I}_1 \rightarrow
{\mathcal I}_2$ of sheaves of ${\mathcal O}_X$-modules extends to ${\mathcal
I}_1 \otimes {\mathcal K}_X \rightarrow {\mathcal I}_2 \otimes {\mathcal
K}_X$.  Each of these is isomorphic to ${\mathcal K}_X$, so the map is given by
multiplication by a global section $f \in {\mathcal K}_X$.  If ${\mathcal F}$
is coherent satisfying $S_2$ and locally free of rank $1$ at every generic
point, then ${\mathcal F} \otimes {\mathcal K}_X \cong {\mathcal K}_X$ and the
natural map ${\mathcal F} \rightarrow {\mathcal F} \otimes {\mathcal K}_X$
makes ${\mathcal F}$ into a nondegenerate fractional ideal.

\bigskip
\noindent
{\bf Warning 2.5.} The usual theory of the sheaf ${\mathcal L}(D) = {\mathcal
I}^{-1}$ associated to a divisor $D$ \cite[2.8]{GD} does not extend to divisors
that may not be reflexive.  However, we can get an analogue of \cite[2.10]{GD} using
the sheaf ${\mathcal M}(D) = {\mathcal H}om({\mathcal I},\omega)$.  Note that ${\mathcal M}(D)$ is $\omega$-reflexive by $(1.6)$ and therefore satisfies $S_2$.

\bigskip
\noindent
{\bf Proposition 2.6.} {\em Let $X$ be a Cohen--Macaulay scheme with canonical
sheaf $\omega$, and for any divisor $D$, corresponding to an ideal sheaf
${\mathcal I}$, let ${\mathcal M}(D) = {\mathcal H}om({\mathcal I},\omega)$. 
If $D$ is an effective divisor, denoting also by $D$ the associated closed
subscheme, there are two natural exact sequences}
\[
0 \rightarrow {\mathcal I} \rightarrow {\mathcal O}_X \rightarrow {\mathcal
O}_D \rightarrow 0
\]
{\em and}
\[
0 \rightarrow \omega_X \rightarrow {\mathcal M}(D) \rightarrow \omega_D
\rightarrow 0.
\]

\bigskip
\noindent
{\bf {\em Proof.}} The first is the defining sequence of $D$.  The second is obtained
by applying ${\mathcal H}om(\cdot,\omega_X)$ to the first and noting (since $X$
is Cohen--Macaulay) that $\omega_D \cong {\mathcal E}xt_{{\mathcal
O}_X}^1({\mathcal O}_D,\omega_X)$.

\bigskip
\noindent
{\bf Definition--Remark 2.7.} Even though $X$ has a canonical sheaf $\omega_X$, it may not
have a canonical divisor.  By {\em canonical divisor} we mean a generalized divisor $K$ whose ideal satisfies ${\mathcal I}_K^{-1} \cong \omega_X$.  Since the ideal of any divisor is locally free at the generic points, the existence of a canonical divisor implies that $X$ satisfies $G_0$.  In this case we see also that $\omega_X$ must be reflexive, by \cite[1.8]{GD}.  Since there are schemes satisfying $G_0$ and $S_2$ on which $\omega$ is not reflexive $(1.9)$, we conclude that $G_0$ and $S_2$ are not sufficient conditions for the existence of a canonical divisor.

However, if 
 $X$ satisfies $G_0$ and $S_2$, then $\omega_X$ is
locally free of rank $1$ at every generic point, so is isomorphic to a fractional
ideal.  We choose and fix an embedding $\omega_X \subseteq {\mathcal K}_X$, and
call the corresponding divisor $M_X$ the {\em anticanonical} divisor.  As a
divisor it depends on the choice of embedding $\omega_X \subseteq {\mathcal
K}_X$, but is unique up to linear equivalence.  If $X$ satisfies in addition
$G_1$, then $\omega$ is invertible in codimension $1$, so we can define a
{\em canonical} divisor $K = -M$, which will be an almost Cartier
divisor.

\bigskip
\noindent
{\bf Definition.} For any two divisors $D_1,D_2$, we define $D_1(-D_2)$ to be
the divisor whose sheaf of ideals is ${\mathcal H}om({\mathcal I}_2,{\mathcal
I}_1)^{\sim}$.  In general, this operation may not be well-behaved, but we do
have the following.

\bigskip
\noindent
{\bf Proposition 2.8.} {\em The operation $D_1(-D_2)$ has the following
properties.}
\begin{itemize}
\item[(a)] {\em $0(-D_2) = -D_2$ and $D_1(-0) = D_1$.}
\item[(b)] {\em If $E$ is almost Cartier, $(D_1+E)(-D_2) = D_1(-(D_2-E)) =
D_1(-D_2) + E$.}
\item[(c)] {\em In particular, if either $D_1$ or $D_2$ is almost Cartier, then
$D_1(-D_2) = D_1 + (-D_2)$.}
\item[(d)] {\em If $X$ satisfies $G_0$, and $D_1 \sim M + E$, where $M$ is the
anticanonical divisor and $E$ is almost Cartier, then $D_1(-D_1(-D_2)) = D_2$
for any $D_2$.}
\end{itemize}

\bigskip
\noindent
{\bf {\em Proof.}} (a), (b), (c) are immediate, since an almost Cartier divisor is
invertible in codimension $1$, and equality of divisors can be tested in
codimension $1$.  (d) corresponds to the fact that any divisor has an ideal
sheaf satisfying $S_2$, and hence is $\omega$-reflexive $(1.5)$.

\bigskip
\noindent
{\bf Remark 2.9.}  We take this opportunity to point out an error in \cite[2.9]{GD}.  Assuming that $X$ satisfies $G_1$ and $S_2$ as in that paper, it is true that every nondegenerate section $s \in \Gamma(X,{\mathcal L}(D))$ gives rise to an effective divisor $D'$ in the complete linear system $|D|$, and all $D'$ arise in this way.  Two sections $s_1$ and $s_2$ give rise to the same divisor $D'$ if and only if they differ by an isomorphism of ${\mathcal L}(D)$.  If $D$ is almost Cartier, the isomorphisms of ${\mathcal L}(D)$ are given by sections of $\Gamma(X,{\mathcal O}_X^*)$ as stated there.  So in the familiar case of $X$ integral projective, $\Gamma(X,{\mathcal O}_X^*) = k^*$ and $|D|$ is simply the projective space associated to the vector space $\Gamma(X,{\mathcal L}(D))$.

Suppose, however, that $D$ is not almost Cartier.  Then there may be more isomorphisms of ${\mathcal L}(D)$ and the statement of \cite[2.9]{GD} is not correct.  For example, let $X = L_1 \cup L_2$ be the union of two lines in ${\mathbb P}^2$ meeting at a point $P$, and let $D$ be the divisor $P$.  Then one can verify that $\dim \Gamma(X,{\mathcal L}(D)) = 2$, and $\mbox{Isom}({\mathcal L},{\mathcal L}) = k^* \oplus k^*$, so that the complete linear system $|D|$ consists just of the single divisor $D$, as we expect.  (Cf.~\cite[3.3]{GD} for a relevant calculation.)

How does this discussion extend to the case of the present paper, where $X$ is only assumed to satisfy $S_2$?  We cannot use the sheaf ${\mathcal L}(D)$.  Instead, for each effective divisor $D' \sim D$, we take $\omega$-duals of ${\mathcal I}_{D'} \subseteq {\mathcal O}_X$ to get $\omega_X \subseteq {\mathcal H}om({\mathcal I}_{D'},\omega_X) \cong {\mathcal M}(D)$, and this gives a section $s$ of the sheaf ${\mathcal N}(D) = {\mathcal H}om(\omega_X,{\mathcal M}(D))$.  Conversely, nondegenerate sections of ${\mathcal N}(D)$ give effective divisors $D' \sim D$ by reversing the process.  The ambiguity of $s$ is again in $\mbox{Isom}({\mathcal N}(D),{\mathcal N}(D)) \cong \mbox{Isom}({\mathcal I}_D,{\mathcal I}_D)$.

\section{Biliaison}
\label{sec3}

In this section we generalize the notion of biliaison introduced in \cite[\S
4]{GD} and \cite[\S 5.4]{M}.  Note that the word biliaison is not a synonym for
even liaison.  We also generalize the results of \cite[\S 5]{KMMNP} so as to
remove the $G_1$ hypotheses.  In fact, it was the attempt to put those results
in a more natural context that led to this paper.

\bigskip
\noindent
{\bf Definition.} Let $V_1$ and $V_2$ be equidimensional closed subschemes of
dimension $r$ of ${\mathbb P}_k^n$.  We say that $V_2$ is obtained by an {\em
elementary biliaison} of height $h$ from $V_1$ if there exists an ACM scheme
$X$ in ${\mathbb P}^n$, of dimension $r+1$, containing $V_1$ and $V_2$, and so
that $V_2 \sim V_1 + hH$ as generalized divisors on $X$, where $H$ denotes the
hyperplane class.  The equivalence relation generated by elementary biliaisons
will be called {\em biliaison}.

If we restrict the schemes $X$ in the definition all to be complete
intersection schemes, we will speak of CI-{\em biliaison}.  If we restrict the
schemes $X$ to be ACM schemes satisfying $G_0$, will speak of {\em Gorenstein
biliaison} or G-{\em biliaison}.

\bigskip
\noindent
{\bf Remark 3.1.} As was shown in \cite[4.4]{GD} the relation of CI-biliaison
is equivalent to even CI-liaison in the usual sense.

\bigskip
\noindent
{\bf Proposition 3.2.} {\em Suppose $V_2$ is obtained from $V_1$ by an
elementary biliaison of height $h$ on $X$, with $\dim V_1 = \dim V_2 = r$.}
\begin{itemize}
\item[a)] {\em Then reciprocally, $V_1$ is obtained from $V_2$ by an elementary
biliaison of height $-h$.}
\item[b)] {\em The higher Rao modules $M_V^i = H_*^i({\mathcal I}_{V,{\mathbb
P}^n})$ are related as follows:
\[
M_{V_2}^i \cong M_{V_1}^i(-h) \mbox{ for } 1 \le i \le r.
\]
}
\item[c)] {\em The Hilbert polynomials are related by}
\[
\chi({\mathcal O}_{V_2}(m)) = \chi({\mathcal O}_X(m)) - \chi({\mathcal
O}_X(m-h)) + \chi({\mathcal O}_{V_1}(m-h)).
\]
\end{itemize}

\bigskip
\noindent
{\bf {\em Proof.}} a)  If $V_2 \sim V_1 + hH$ then $V_1 \sim V_2 - hH$.

b) and c) have the same proof as \cite[4.5]{GD} since only the ACM property of
$X$ was used there.

\bigskip
\noindent
{\bf Lemma 3.3.} {\em Let $X$ be an {\em ACM} scheme satisfying $G_0$ in
${\mathbb P}^n$.  Let $Y \subseteq X$ be an effective divisor, $Y \sim M + mH$,
where $M$ is the anticanonical divisor and $H$ is the hyperplane divisor.  Then
$Y$ is an arithmetically Gorenstein {\em (AG)} scheme in ${\mathbb P}^n$.}

\bigskip
\noindent
{\bf {\em Proof.}} Let $X$ be of dimension $r+1$ so that $Y$ is of dimension $r$.  To
show that $Y$ is AG is equivalent to showing that $Y$ is ACM and $\omega_Y
\cong {\mathcal O}_Y(\ell)$ for some $\ell \in {\mathbb Z}$.

First to show $Y$ is ACM, we must show $H_*^i({\mathcal I}_{Y,{\mathbb P}^n}) =
0$ for $1 \le i \le r$.  From the exact sequence
\[
0 \rightarrow {\mathcal I}_{X,{\mathbb P}^n} \rightarrow {\mathcal
I}_{Y,{\mathbb P}^n} \rightarrow {\mathcal I}_{Y,X} \rightarrow 0
\]
and the fact that $X$ is ACM, so that $H_*^i({\mathcal I}_{X,{\mathbb P}^n}) =
0$ for $1 \le i \le r+1$, it is equivalent to show $H_*^i({\mathcal
I}_{Y,X})  = 0$ for $1 \le i \le r$.  Now $Y \sim M + mH$ by hypothesis, so
${\mathcal I}_{Y,X} \cong \omega_X(-m)$.  By Serre duality on $X$,
$H_*^i(\omega_X(-m))$ is dual to $H_*^{r+1-i}({\mathcal O}_X(m))$.  These
latter are $0$ for $1 \le i \le r$ since $X$ is ACM.  Hence $Y$ is ACM.

To study the canonical sheaf $\omega_Y$, we use the second exact sequence of
$(2.6)$, namely
\[
0 \rightarrow \omega_X \rightarrow {\mathcal M}(Y) \rightarrow \omega_Y
\rightarrow 0.
\]
Now since ${\mathcal I}_Y \cong \omega_X(-m)$, we have $\omega_X \cong
{\mathcal I}_Y(m)$ and ${\mathcal M}(Y) = {\mathcal H}om({\mathcal
I}_Y,\omega_X) \cong {\mathcal O}_X(m)$.  Thus $\omega_Y \cong {\mathcal
O}_Y(m)$ and $Y$ is arithmetically Gorenstein.

\bigskip
\noindent
{\bf Remark 3.4.} An algebraic version of this result was given in
\cite[5.2]{KMMNP}, and a geometric version with the added hypothesis $G_1$ in
\cite[5.4]{KMMNP}.

\bigskip
\noindent
{\bf Definition.}  Two subschemes $V_1$ and $V_2$ of ${\mathbb P}^n$, equidimensional of the same dimension and without embedded components are {\em linked} by a scheme $Y$ if $Y$ contains $V_1$ and $V_2$ and ${\mathcal I}_{V_i,Y} \cong {\mathcal H}om({\mathcal O}_{V_j},{\mathcal O}_Y)$ for $i,j = 1,2,i \ne j$.  If $Y$ is a complete intersection, it is called a CI-{\em linkage}; if $Y$ is arithmetically Gorenstein, it is a {\em Gorenstein linkage}.  If $Y$ is an arithmetically Gorenstein scheme of the form $M+mH$ on some ACM scheme $X$ satisfying $G_0$ (as in $(3.3)$ above), then we will say it is a {\em strict Gorenstein linkage}.  (This is a slight generalization of the terminology of \cite[\S1]{H3}, where we required that $X$ should satisfy $G_1$.)

The equivalence relation generated by CI-linkages is CI-{\em liaison}, by Gorenstein linkages, {\em Gorenstein liaison}, and by strict Gorenstein linkages, {\em strict Gorenstein liaison}.  If the liaison can be accomplished by an even number of linkages, then it is an {\em even} CI-liaison (resp.\ Gorenstein liaison, resp.\ strict Gorenstein liaison).

\bigskip
\noindent
{\bf Theorem 3.5.} {\em Suppose that $V_2$ is obtained from $V_1$ by an
elementary biliaison on an {\em ACM} scheme $X$ satisfying $G_0$.  Then $V_2$
can be obtained from $V_1$ by two strict Gorenstein linkages.}

\bigskip
\noindent
{\bf {\em Proof.}} The proof is almost the same as \cite[4.3]{GD}, transposed
into our context.  We assume that $V_2 \sim V_1 + hH$.  Thus there is a
principal divisor $(f)$ such that $V_2 = V_1 + hH + (f)$.  Taking $M$ to be the
anticanonical divisor and using $(2.8)$, we can write
\[
M(-V_1) = M(-V_2) + hH + (f).
\]
Now by \cite[2.11]{GD}, which still holds in our case, we can find an effective
Cartier divisor $E \sim mH$ such that $W = E + M(-V_1) = (M+E)(-V_1)$ is
effective.  Now let $Y = M+E$.  Then $Y$ is an effective divisor that is
arithmetically Gorenstein by $(3.4)$, and I claim that $V_1$ and $W$ are
Gorenstein linked by $Y$.  Indeed, the same argument as in the proof of
\cite[4.1]{GD} shows that ${\mathcal I}_{W,Y} \cong {\mathcal H}om({\mathcal
O}_{V_1},{\mathcal O}_Y)$.  Since by $(2.8)$ we also have $V_1 = (M+E)(-W)$, we
obtain the reverse isomorphism ${\mathcal I}_{V_1,Y} \cong {\mathcal
H}om({\mathcal O}_W,{\mathcal O}_Y)$, so $V_1$ and $W$ are linked by $Y$.

We also have $W = E+M(-V_2)+hH+(f)$, so if we let $Y' = E+M+hH+(f)$, then $Y'$
will also be an effective divisor that is an arithmetically Gorenstein scheme,
and as above, we see that $W$ and $V_2$ are linked by $Y'$.  Thus $V_2$ is
obtained from $V_1$ by two strict Gorenstein linkages.

\bigskip
\noindent
{\bf Corollary 3.6.} {\em Every Gorenstein biliaison is an even strict Gorenstein
liaison.}

\bigskip
\noindent
{\bf Remark 3.7.} This theorem was proved for a trivial elementary biliaison
$V_2 = V_1 + hH$ (with no linear equivalence) in \cite[5.10]{KMMNP}, and with
the extra hypothesis $G_1$ in \cite[5.14]{KMMNP}.

\bigskip
\noindent
{\bf Remark 3.8.} In section~\ref{sec5} below we will prove a converse to this theorem in codimension~$3$.

\bigskip
\noindent
{\bf Example 3.9.} Let $P$ be a point in ${\mathbb P}^3$, and let $X$ be the
union of three non-coplanar lines through $P$.  Then $X$ satisfies $G_0$ but
not $G_1$.  If $H$ is a hyperplane section of $X$ containing $P$, then $V = P +
H$ is the divisor defined by the square of the ideal of $P$.  Thus $P$ and $V$
are related by one $G$-biliaison, and hence are evenly $G$ linked.  Cf.\
\cite[4.1]{KMMNP}, where this was proved by a different method.

\section{The theorem of Gaeta}
\label{sec4}

To illustrate the theory of biliaison, we give a new proof the theorem of KMMNP--Gaeta
\cite[3.6]{KMMNP}.  The statement given there is that every standard
determinantal scheme is glicci.  We prove a slightly stronger result.

\bigskip
\noindent
{\bf Theorem 4.1.} {\em Every standard determinantal scheme in ${\mathbb P}^n$
can be obtained from a linear variety by a finite number of ascending Gorenstein
biliaisons.  In particular, it is glicci by $(3.3)$.}

\bigskip
\noindent
{\bf {\em Proof.}}  We follow the terminology and notation of \cite[3.6]{KMMNP}.  Let
$V \subseteq {\mathbb P}^n$ be a standard determinantal scheme, i.e., a scheme of codimension $c+1$ whose ideal
$I_V$ is generated by the $t \times t$ minors of a $t \times (t+c)$ homogeneous
matrix $A$ for some $t > 0$.  Let $B$ be the matrix obtained by omitting the last column of
$A$.  Then $V$ is contained in the determinantal scheme $S$ defined by the $t
\times t$ minors of $B$.  By Step~I of the proof of \cite[3.6]{KMMNP}, $S$ is
good determinantal.  Hence it is generically a complete intersection
\cite[3.2]{KMMNP}, and so satisfies $G_0$.

Let $A'$ be the matrix obtained by omitting the last row of $B$.  Then $V'$,
defined by the $(t-1) \times (t-1)$ minors of $A'$, is also contained in $S$.  We
will show that $V \sim V' + mH$ on $S$ for some $m > 0$, so that $V$ is obtained
by an ascending elementary Gorenstein biliaison from $V'$.  Continuing in this
manner, after a finite number of $G$-biliaisons, we reduce to the case $t=1$,
when $V$ is a complete intersection.  From these one can perform descending
$CI$-biliaisons to a linear variety.

Let $R$ be the homogeneous coordinate ring of ${\mathbb P}^n$, and let $R_S =
R/I_S$ be the homogeneous coordinate ring of $S$.  The ideal of $V$ in $S$ is
generated by the images in $R_S$ of the $t \times t$ minors of $A$ that include
the last column.  The $t \times t$ minors that do not include the last column
are just the generators of $I_S$.  On the other hand, the ideal of $V'$ in $S$
is generated by the images of the $(t-1) \times (t-1)$ minors of $A'$.  So there
is a one-to-one correspondence between generators $N$ of $V$ in $S$ and
generators $N'$ of $V'$ in $S$, obtained by omitting the last row and column of
the corresponding $t \times t$ matrix.  We will show that the quotient $N/N'$ of
corresponding generators is an element of $H^0({\mathcal K}_S(m))$, independent
of the choice of
$N$, where
${\mathcal K}_S$ is the sheaf of total quotient rings of $S$, and $m$ is the
difference in degrees of $N$ and $N'$.   This will show that ${\mathcal I}_{V,S}
\cong {\mathcal I}_{V',S}(-m)$, and so we have the desired biliaison.  Note that $m$ is the degree of the element in the lower right-hand corner of the original matrix $A$.

To show that $N/N'$ is independent of the choice of $N(\mbox{mod } I_S)$, it
will be sufficient to compare two such that differ by one column only.  So let
$M$ be a $t \times t$ minor of $B$, let $N_1$ be obtained by deleting the first
column of $M$ and adding the last column of $A$; let $N_2$ be obtained by
deleting the second column of $M$ and adding the last column of $A$.  Then $N_1$
and $N_2$ are two generators of $I_V$, and the corresponding generators
$N'_1,N'_2$ of $I_{V'}$ are just $M_{t1}$ and $M_{t2}$, where $M_{ij}$ denotes
the minor of $M$ obtained by deleting the $i^{th}$ row and the $j^{th}$ column. 
We need to show that $N_1/N'_1 = N_2/N'_2 \mbox{ mod } I_S$.  By making general
row and column operations on $A$ at the beginning, we may assume that all the
$N'_i$ are non-zero-divisors in $R_S$.  So we must show that $N_1N'_2 - N_2N'_1
\in I_S$.  

Let the last column of $A$ be $u_1,\dots,u_t$.  We will expand $N_1$ and $N_2$
along this last column.  The coefficient of $u_t$ in $N_1N'_2 - N_2N'_1$ is just
$N'_1N'_2 - N'_2N'_1 = 0$.  For $i \ne t$, the coefficient of $u_i$ is
$M_{i1}M_{t2} - M_{i2}M_{t1}$.  The proof is then completed by the following
identity among determinants, since $M \in I_S$.

\bigskip
\noindent
{\bf Lemma 4.2.} {\em Let $M$ be a $t \times t$ matrix, let $M_{ij}$ denote the
minor obtained by deleting the $i^{th}$ row and the $j^{th}$ column; let
$M_{ik,jl}$ denote the minor obtained by deleting the $i^{th}$ and $k^{th}$ rows
and the $j^{th}$ and $l^{th}$ columns.  Then the determinants satisfy}
\[
M_{ij} \cdot M_{kl} - M_{il} \cdot M_{kj} = \pm M_{ij,kl} \cdot M.
\]

\bigskip
\noindent
{\bf {\em Proof.}}  \cite[p.~132ff]{Muir}.

\bigskip
\noindent
{\bf Example 4.3.}  The $4 \times 4$ minors of a general $4 \times 6$ matrix of
linear forms in ${\mathbb P}^4$ define an irreducible smooth curve $C$ of degree
$20$ and genus $26$ which, according to the theorem, can be obtained by
ascending Gorenstein biliaisons from a line.  However these curves are not
general in the Hilbert scheme, and it is known that a general smooth curve of
degree $20$ and genus $26$ is ACM, but cannot be obtained by ascending
biliaisons from a line.  It is unknown whether it is glicci \cite[3.9]{H3}.

\section{Strict Gorenstein liaison}
\label{sec5}

The main result of this section is a converse to $(3.6)$ in codimension~$3$.

\bigskip
\noindent
{\bf Theorem 5.1.} {\em For subschemes of codimension~$3$ in ${\mathbb P}^n$ (equidimensional and without embedded components), any even strict Gorenstein liaison is a Gorenstein biliaison.}

\bigskip
\noindent
{\bf {\em Proof.}}  Suppose $V$ and $V'$ of codimension~$3$ in ${\mathbb P}^n$ are related by even strict Gorenstein liaison.  Then there is a sequence
\[
V = V_0,V_1,V_2,\dots,V_{2k} = V'
\]
for some $k$, where each $V_i$ is related to $V_{i+1}$ by a strict Gorenstein linkage.  By composition of biliaisons, it will be sufficient to treat the case $k=1$, i.e., when there is just one intermediary scheme $Z$, and $V$ to $Z$ is a strict Gorenstein linkage by $Y$ of the form $M+mH$ on an ACM scheme $X$ satisfying $G_0$, and $Z$ to $V'$ similarly is linked by a $Y'$ of the form $M' + m'H'$ on an ACM scheme $X'$ satisfying $G_0$.

Since $X$ and $X'$ are both ACM of codimension~$2$ in ${\mathbb P}^n$, they are in the same CI-biliaison class, by the classical Gaeta's theorem \cite[6.1.4]{M}.  Thus we can apply Lemma~$5.2$ (below) and find a chain
\[
X = X_0,X_1,\dots,X_r = X'
\]
of ACM schemes satisfying $G_0$ and each containing $Z$, such that each $X_i$ is directly CI-linked to $X_{i+1}$, and $X_i$ and $X_{i+1}$ have no common components.

Now for each $i = 1,\dots,r$, let $D_i = X_{i-1} \cap X_i$.  By Lemma~$5.3$ (below), $D_i$ is an AG scheme of the form $M+mH$ on $X_{i-1}$ and on $X_i$.  Since the $X_i$ all contain $Z$, so do the $D_i$.  For each $i = 1,\dots,r$, let $W_i$ be the scheme linked to $Z$ by $D_i$.  We consider the chain of strict Gorenstein linkages
\[
V = W_0,Z,W_1,Z,W_2,Z,\dots,W_r,Z,V' = W_{r+1}.
\]
Here, for each $i = 0,\dots,r$, the two links $W_i,Z,W_{i+1}$ are both strict Gorenstein links on the same ACM scheme $X_i$.  Now, as in the proof of \cite[4.1]{GD} we see that $W_i$ being linked to $Z$ by $M+mH$ on $X_i$ is equivalent to saying $Z \sim (M+mH)(-W_i)$ on $X_i$.  Similarly, $Z$ linked to $W_{i+1}$ by $M+m'H$ on $X_i$ says $W_{i+1} \sim (M+m'H)(-Z)$.  Substituting the first expression for $Z$ in the second expression for $W_{i+1}$, we find using $(2.8)$ that $W_{i+1} \sim W_i + (m'-m)H$ on $X_i$, which is a single Gorenstein biliaison.

Thus $V$ is joined to $V'$ by the chain of Gorenstein biliaisons
\[
V = W_0,W_1,\dots,W_r,W_{r+1} = V'.
\]

\bigskip
\noindent
{\bf Lemma 5.2.} {\em Suppose given $X,X'$ locally Cohen--Macaulay subschemes of codimension~$2$ in ${\mathbb P}^n$, both satisfying $G_0$, and both containing a given closed subscheme $Z$ of codimension at least $3$ in ${\mathbb P}^n$, and with $X,X'$ in the same {\em CI}-liaison class in ${\mathbb P}^n$.  Then there exists a chain
\[
X = X_0,X_1,\dots,X_r = X'
\]
of locally Cohen--Macaulay subschemes of ${\mathbb P}^n$, each containing $Z$, such that each $X_{i+1}$ is obtained by a single geometric {\em CI}-linkage from $X_i$.  In particular, $X_i$ and $X_{i+1}$ will have no common components, so each will be generically locally complete intersection and therefore will satisfy $G_0$.}

\bigskip
\noindent
{\bf {\em Proof.}}  Note first that the hypothesis $G_0$ implies that $X$ and $X'$ are generically locally complete intersection, since they are in codimension~$2$.  If $X$ to $X'$ is an odd liaison, we can  make a single general geometric liaison from $X'$ to a new $X''$ also containing $Z$, and thus reduce to the case of an even liaison.  Then, since $X$ and $X'$ are in the same even liaison class, by Rao's theorem, they have ${\mathcal N}$-type resolutions with stably equivalent sheaves ${\mathcal N}_i$ up to twist.  By adding dissoci\'e sheaves, we can write ${\mathcal N}$-type resolutions
\[
\begin{array}{cccccccccc}
0 &\to &{\mathcal L} &\to &{\mathcal N} &\to &{\mathcal I}_X(a) &\to &0 \\
0 &\to &{\mathcal L}' &\to &{\mathcal N} &\to &{\mathcal I}_{X'}(a') &\to &0
\end{array}
\]
with the same locally free sheaf ${\mathcal N}$ in the middle, and ${\mathcal L},{\mathcal L}'$ dissoci\'e.

Now we will follow the plan of the proof of \cite[3.1]{CH} to obtain a chain
\[
X = X_0,X_2,X_4,\dots,X_{2k} = X'
\]
of locally Cohen--Macaulay subschemes containing $Z$, such that for each $i$, $X_{2i}$ and $X_{2i+2}$ have no common components and are related by a single elementary CI-biliaison on a hypersurface $S_i$.

Write ${\mathcal L} = \oplus {\mathcal L}_i$ with ${\mathcal L}_i$ invertible, $i = 1,\dots,t$.  Since $X$ is generically locally complete intersection, the rank of the map
\[
{\mathcal L}(\xi) \to {\mathcal N}(\xi)
\]
is $t-1$ for each generic point $\xi$ of $X$.  Thus, reordering if necessary, we define ${\mathcal F}$ by
\[
0 \to \bigoplus_{i \ge 2} {\mathcal L}_i \to {\mathcal N} \to {\mathcal F} \to 0
\]
and ${\mathcal F}$ will be torsion-free of rank~$2$, and locally free at each generic point $\xi$ of $X$.  Now choose $b \gg 0$ so that ${\mathcal I}_Z \otimes {\mathcal N}(b)$ is generated by global sections and take $s_1 \in H^0({\mathcal I}_Z \otimes {\mathcal N}(b))$ a sufficiently general section.  Let $Y_1$ be defined by
\[
0 \to {\mathcal O}(-b)\ {\overset {s_1}{\to}}\ {\mathcal F} \to {\mathcal I}_{Y_1}(a_1) \to 0.
\]
Then $Y_1$ contains $Z$, and $Y_1$ has no component in common with $X$, and $Y_1$ is obtained from $X$ by a single CI-biliaison \cite[3.3]{CH}.  Furthermore, we can lift $s_1$ to ${\mathcal N}$ in such a way that $s_1(\xi_1) \ne 0$ in ${\mathcal N}(\xi_1)$ for each generic point of $Y_1$.  In terms of ${\mathcal N}$ we now have
\[
0 \to {\mathcal O}(-b) \oplus \bigoplus_{i \ge 2} {\mathcal L}_i \to {\mathcal N} \to {\mathcal I}_{Y_1}(a_1) \to 0.
\]

We repeat this process with each ${\mathcal L}_i$ in turn, obtaining a sequence of biliaisons $X,Y_1,Y_2,\dots,Y_t$, each one containing $Z$ and having no components in common with its neighbors.

We do the same thing with $X'$, obtaining a similar sequence $X',Y'_1,\dots,Y'_t$.  Then we observe that one can take the same $b$ in both cases, and since the sections $s_1,\dots,s_t,s'_1,\dots,s'_t$ are all sufficiently general, we can take $s_i = s'_i$ for each $i$, and thus $Y_t = Y'_t$.  This connects $X$ and $X'$ by biliaisons, all containing $Z$.  Now just relabel $Y_i$ and $Y'_i$ as $X_{2j}$ to get the sequence of biliaison above.

To conclude, let $X_2$ and $X_4$ for example be a biliaison on a hypersurface $S$, where $X_2,X_4$ both contain $Z$ and have no common component.  Then $X_2$ and $X_4$ are both generically Cartier divisors on $S$.  When we link them both to a divisor $X_3$ on $S$, as in the proof of $(3.5)$, we can take $X_3$ to have no component in common with $X_2$ and $X_4$, and by adding a complete intersection on $S$ containing $Z$ if necessary, we may assume $X_3$ contains $Z$.  Thus the sequence of biliaisons connecting $X$ and $X'$ can be filled in to a sequence of geometric liaisons as required.

\bigskip
\noindent
{\bf Lemma 5.3.} {\em Let $X_1,X_2$ be {\em ACM} schemes in ${\mathbb P}^n$ that have no common component and are directly linked by an {\em AG} scheme $S$.  Then $D = X_1 \cap X_2$ is arithmetically Gorenstein; moreover, it is of the form $M+\ell H$ on each of $X_1,X_2$, where $\ell$ is the integer for which $\omega_S \cong {\mathcal O}_S(\ell)$.}

\bigskip
\noindent
{\bf {\em Proof}} (cf.\ \cite[4.2.1]{M}).  The fact that $D$ is ACM follows from the exact sequence
\[
0 \to {\mathcal O}_S \to {\mathcal O}_{X_1} \oplus {\mathcal O}_{X_2} \to {\mathcal O}_D \to 0.
\]
Since $S$ is AG, its dualizing sheaf $\omega_S$ is isomorphic to ${\mathcal O}_S(\ell)$ for some $\ell \in {\mathbb Z}$.  Because of the linkage, ${\mathcal I}_{X_1,S} \cong {\mathcal H}om({\mathcal O}_{X_2},{\mathcal O}_S)$.  Note that ${\mathcal I}_{X_1,S} = {\mathcal I}_{D,X_2}$ by a standard isomorphism theorem for ideals.  Since $\omega_{X_2} = {\mathcal H}om({\mathcal O}_{X_2},\omega_S)$, we find that ${\mathcal I}_{D,X_2} \cong \omega_{X_2}(-\ell)$.  This says $D \sim M + \ell H$ on $X_2$.  The same argument shows that $D \sim M + \ell H$ on $X_1$ also.

\section{Conclusion}
\label{sec6}

If we reflect on the outstanding problem whether every ACM subscheme of ${\mathbb P}^n$ is glicci, we can appreciate the usefulness of the extended notion of generalized divisors in this paper.  It has allowed us to prove the theorem of KMMNP--Gaeta in a strengthened form, namely that any standard determinantal scheme in ${\mathbb P}^n$ can be obtained by ascending Gorenstein biliaisons from a linear space.  This also makes clear the special nature of determinantal schemes, since there are known examples of other ACM schemes that cannot be obtained by ascending Gorenstein biliaisons from a linear space, even though it is still unknown whether they are glicci or not (for curves in ${\mathbb P}^4$, see \cite[3.9]{H3}, and for points in ${\mathbb P}^3$ see \cite[7.2]{HSS}).

We also observe that in most known proofs that some class of ACM schemes is glicci (such as the theorem of KMMNP--Gaeta discussed here) the proof could be accomplished using Gorenstein biliaisons, hence using only strict Gorenstein liaisons.  Since there are AG schemes in ${\mathbb P}^n$ not of the special form $M + mH$ on some ACM scheme of one dimension higher (for curves in ${\mathbb P}^4$, see \cite[3.6, 3.11]{GAGC} and for points in ${\mathbb P}^3$ see \cite[3.4, 6.8]{HSS}), this suggests that it would be worthwhile to investigate more deeply what kind of $G$-liaisons can be accomplished using AG schemes not of this special form.

\end{document}